\documentclass[12pt]{amsart}
\usepackage{amsmath,amsthm,amsfonts,amssymb}
\usepackage[dvips, usenames, dvipsnames]{xcolor}
\definecolor{orange2}{HTML}{FF7F2A}
\usepackage{pinlabel}

\usepackage[colorlinks=true,linkcolor=blue,citecolor=blue,urlcolor=blue]{hyperref}

\usepackage{tikz}
\usetikzlibrary{arrows,automata}
\usepackage[latin1]{inputenc}
\tikzset{cdlabel/.style={above,sloped,
    execute at begin node=$\scriptstyle,execute at end node=$}}    
\tikzset{algarrow/.style={->, thick}}   
\def\mathcenter#1{%
  \vcenter{\hbox{$#1$}}%
}

\normalfont\upshape

\usepackage[lmargin=1.2in,rmargin=1.2in,tmargin=1.2in,bmargin=1.2in]{geometry}

\newtheorem{theorem}{Theorem}[section]

\newtheorem{definition}[theorem]{Definition}

\newtheorem*{remark*}{Remark}
\newtheorem*{example*}{Example}

\newcommand{\R}{\ensuremath{\mathbb{R}}}
\newcommand{\Z}{\ensuremath{\mathbb{Z}}}

\newcommand\HFK{{\rm {HFK}}}
\newcommand\HFKh{{\rm {\widehat{HFK}}}}
\newcommand\HFKm{{\rm {HFK^{-}}}}

\newcommand\CFKh{{\rm {\widehat{CFK}}}}
\newcommand\CFKm{{\rm {CFK^{-}}}}

\def\x{\mathbf{x}}
\def\y{\mathbf{y}}

\def\gen{\mathfrak{S}}

\newcommand{\ctt}{\widetilde{\rm {CT}}}
\newcommand{\ctm}{{\rm {CT}}^-}

\newcommand{\bdy}{\partial}

\def\Ft{\mathbb{F}_2} 
\newcommand{\eee}[1]{e_{#1}} 

\newcommand{\eeA}[1]{\eee{#1}^A}
\newcommand{\eeD}[1]{\eee{#1}^D}
\newcommand{\T}{\mathcal{T}} 
\newcommand{\Tdec}{\mathbb{T}} 

\newcommand{\tDA}{\mathit{DA}}

\newcommand{\sss}{{\bf s}}
\newcommand{\ttt}{{\bf t}}

\newcommand{\HH}{\mathcal{H}}
\newcommand{\balpha}{\boldsymbol\alpha}
\newcommand{\bbeta}{\boldsymbol\beta}
\newcommand{\zzz}{\mathbf{z}}
\newcommand{\www}{\mathbf{w}}

\newcommand{\ophat}[1]{\widehat{#1}}
\newcommand{\opminus}[1]{#1^-}
\newcommand{\alg}[1]{\mathcal{#1}}
\newcommand{\am}[1]{\opminus{\alg{A}} (#1)}
\newcommand{\ah}[1]{\ophat{\alg{A}} (#1)}
\newcommand{\im}[1]{\opminus{\alg{I}} (#1)}

\begin{document}

\title[An introduction to tangle Floer homology]{An introduction to tangle Floer homology}
\author[PETKOVA and V\'ERTESI]{Ina Petkova, Vera V\'ertesi}

\thanks{}

\address {Department of Mathematics, Columbia University\\ New York, NY 10027}
\email{ina@math.columbia.edu}

\address{Institut de Recherche Math\'ematique Avanc\'ee \\Universit\'e de Strasbourg}
\email{vertesi@math.unistra.fr}

\begin{abstract}
This paper is a short introduction to the combinatorial version of  tangle Floer homology defined in \cite{pv}. There are two equivalent definitions---one in terms of strand diagrams, and one in terms of bordered grid diagrams. We present both, discuss the correspondence, and carry out some explicit computations. 
\end{abstract}

\maketitle


\section{Introduction} 
\label{sec:introduction}

Knot Floer homology is a categorification of the Alexander polynomial. It was introduced by Ozsv\'ath--Szab\'o \cite{hfk} and Rasmussen \cite{jrth} in the early 2000s. One associates a bigraded chain complex $\CFKh(\HH)$ over $\Ft$ to a Heegaard  diagram $\HH = (\Sigma, \balpha, \bbeta, \zzz, \www)$ for a link $L$. The generators are combinatorial and can be read off from the intersections of curves on the Heegaard diagram, whereas the differential counts pseudoholomorphic curves in $\Sigma\times I\times \R$ satisfying certain boundary conditions. The homology of 
$\CFKh(\HH)$ is an invariant of $L$ denoted $\HFKh(L)$. 

Knot Floer homology is a powerful link invariant---it detects genus, detects fiberedness, and an enhanced version called $\HFKm$ contains a concordance invariant $\tau(K)\in\Z$, whose absolute value bounds the 4-ball genus of $K$, and hence the unknotting number of $K$.

Combinatorial versions of knot Floer homology \cite{mos, oszcomb} were  defined soon after the original construction, but they were still global in nature, and our understanding of how local modifications of a knot affect $\HFK$ was very limited \cite{mskein, oszskein}. 

 In \cite{pv}, we ``localize" the construction of knot Floer homology, and define an invariant of oriented tangles. Although  we develop a theory for oriented tangles in general $3$-manifolds with spherical boundaries by using analysis similar to \cite{bfh2, bimod}, in this paper we will focus on a completely combinatorial construction  for tangles in $B^3$ and $I\times S^2$ (we'll think of those as tangles in $I\times \R^2$). 
 
\begin{figure}[h]
\centering
 \includegraphics[scale=0.8]{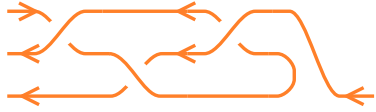} 
       \caption{A projection of a $(3,1)$-tangle $\T$ to $I\times \R$. Here $-\partial^L\T=(-,-,+)$ and $\partial^R\T=(-)$.}\label{fig:tangle}
\end{figure}

An  \emph{$(m,n)$-tangle} $\T$ is a proper, smoothly embedded oriented 1--manifold in  $I\times \R^2$, with boundary $\bdy \T = \bdy^L\T\sqcup \bdy^R\T$, where $\bdy^L\T = \{0\}\times\{1, \ldots, m\}\times \{0\}$ and $\bdy^R\T = \{1\}\times\{1, \ldots, n\}\times \{0\}$, treated as oriented sequences of points; if $m$ or $n$ is zero, the respective boundary is the empty set.  A planar diagram of a tangle is a projection to $I\times\R\times\{0\}$ with no triple intersections, self-tangencies, or cusps, and with over- and under-crossing data preserved (as viewed from the positive $z$ direction).  The boundaries of $\T$ can be thought of as \emph{sign sequences}
$$-\bdy^L\T \in \{+,-\}^m, \quad \bdy^R\T\in\{+,-\}^n,$$
according to the orientation of the tangle at each point ($+$ if the tangle is oriented left-to-right, $-$ if the tangle is oriented right-to-left at that point). See for example Figure \ref{fig:tangle}.
Given two tangles $\T$ and $\T'$ with  $\bdy^R\T=-\bdy^L\T'$, we can concatenate them to obtain a new tangle $\T\circ\T'$, by placing $\T'$ to the right of $\T$.

  We associate a differential graded algebra called $\am{P}$ to a sign sequence $P \in \{+,-\}^n$, and a type $\tDA$ bimodule $\ctm(\Tdec)$ over  $(\am{-\bdy^L\T}$, $\am{\bdy^R\T})$ to a fixed \emph{decomposition} $\Tdec$ of an  $(m,n)$-tangle $\T$. These structures come equipped with  a grading $M$ by $\Z$, called the Maslov grading, and a grading $A$ by $\frac 1 2 \Z$, called the Alexander grading.  Setting certain variables $U_i$ in $\ctm$ to zero, we get a simpler bimodule $\ctt$, which we prove in \cite{pv} to be an invariant of the tangle $\T$ (there is evidence suggesting that $\ctm$ is an invariant too, but we do not at present have a complete proof). Gluing corresponds to taking box tensor product, and for closed links the invariant recovers $\HFK$:
\begin{theorem}\label{thm:main}
Given an $(l,m)$-tangle $\T_1$ with decomposition $\Tdec_1$ and  an $(m,n)$-tangle $\T_2$ with decomposition $\Tdec_2$ with $\bdy^R\T_1 = -\bdy^L\T_2$, let $\Tdec = \Tdec_1 \Tdec_2$ be the corresponding decomposition for the concatenation $\T_1\circ\T_2$. Then there is a bigraded isomorphism
$$\ctm(\Tdec_1)\boxtimes_{\am{\bdy^R\T_1}} \ctm(\Tdec_2)\simeq \ctm(\Tdec).$$
Regarding an $l$-component link $\mathcal L$ (with some decomposition $\mathbb L$) as a $(0,0)$-tangle, we have
$$\ctm(\mathbb L)[-l/2]\{-l/2\}\simeq g\mathit{CFK}^-(\mathcal L)\otimes (\Ft[-1/2]\oplus \Ft[1/2]),$$
where $[i]$ denotes a Maslov grading shift down by $i$, and $\{j\}$ an Alexander grading shift down by $j$.
\end{theorem}
We define $\ctm$ combinatorially, by means of \emph{bordered grid diagrams}, or, equivalently, \emph{strand diagrams}.

\subsection{Outline}

In Section \ref{sec:tangle} we describe the two constructions for $\ctm$ and discuss their correspondence. Section \ref{sec:computations} contains a couple of small explicit computations (a cap and a cup, which glue up to the unknot).

\subsection*{Acknowledgments} The first author thanks the organizers of the 2015 G\"okova Geometry/Topology Conference for an awesome workshop.


\section{Combinatorial tangle Floer homology} 
\label{sec:tangle}


We assume familiarity with the types of algebraic structures discussed in this paper. For some background reading, we suggest \cite[Section 2.1]{pv} (a brief summary) or  \cite[Section 2]{bimod} (a more detailed exposition).

We begin with the definition of the algebra. 

\subsection{The algebra for a sign sequence}\label{ssec:alg}

Let $P = (p_1, \ldots, p_n) \in \{+, -\}^n$ be a sign sequence and let $[n]=\{0,1,\ldots, n\}$. We associate to $P$ a differential graded algebra $\am{P}$ over $\Ft[U_1, \ldots, U_t]$, where the variables $U_1, \ldots, U_t$ correspond to the positively oriented points in $P$. The algebra is generated over $\Ft[U_1, \ldots, U_t]$ by partial bijections $[n] \to [n]$ (i.e. bijections $\sss\to \ttt$ for $\sss, \ttt\subset [n]$), which can be drawn as strand diagrams (up to planar isotopy and Reidemeister III moves), as follows. 

Represent each $p_i$ by a horizontal orange strand  $[0,1]\times \{i-\frac{1}{2}\}$ oriented left-to-right if $p_i=+$  and right-to-left if $p_i = -$ (in \cite{pv}, those are dashed green strands and double orange strands, respectively). Represent a bijection $\phi:\sss\to \ttt$ by black strands connecting $(0,i)$ to $(1, \phi i)$ for $i\in \sss$. We further require that there are no triple intersection points and there are a minimal number of intersection points between strands. 

\begin{figure}[h]
  \centering
  \labellist
  \pinlabel $U_1$ at 522 40
  \endlabellist
  \includegraphics[scale=.6]{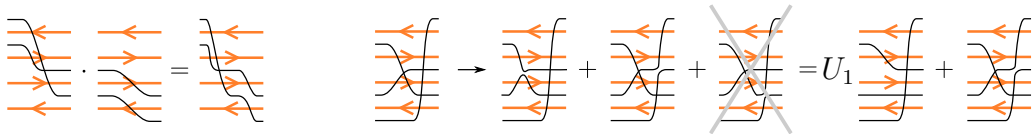}
  \vskip .1cm
  \caption{The algebra $\am{P}$ for $P=(-,+,+,-)$. Left: an example of the multiplication. Right: an example of the differential.}
  \label{fig:alg}
\end{figure}

Let  $a:\sss_1\to \ttt_1, b:\sss_2\to \ttt_2$ be generators. If $\ttt_1 \neq\sss_2$, define the product $ab$ to be $0$. If $\ttt_1 = \sss_2$,  consider the concatenation of a diagram for $a$ to the left and a diagram for $b$ to the right. If there is a black strand that crosses a left-oriented orange strand or another black strand twice, define $ab=0$. Otherwise define $ab = (\prod_i U_i^{n_i}) b\circ a$ where $n_i$ is the number of black strands that ``double cross" the $i^{\mathrm{th}}$ right-oriented orange strand. See Figure \ref{fig:alg}. 

For a generator $a$, define its differential $\bdy a$ as the sum of all ways of locally smoothing one crossing in a diagram for $a$, so that there are no double crossings  between a black strand and a left-oriented orange strand or another black strand, whereas a double crossing between the $i^{\mathrm{th}}$ right-oriented orange strand and a black strand results in a factor of $U_i$, followed by Reidemeister II moves to minimize the total intersection number. See Figure \ref{fig:alg}. 

In other words, product is defined by concatenation, and the differential by smoothing crossings, each followed by  Reidemeister II moves to minimize the total intersection number, where the  relations in Figure \ref{fig:Reid} are applied to the Reidemeister II moves.

The idempotent generators are exactly the identity bijections $\eee{\sss}:\sss\to\sss$, i.e. the strand diagrams consisting of only horizontal strands.  The subalgebra generated over $\Ft[U_1, \ldots, U_t]$ by the idempotents is denoted $\im{P}$.

The algebra $\am{P}$ has a differential grading $M$ called the \emph{Maslov} grading, and an internal grading $A$ called the \emph{Alexander} grading. Those are defined on generators by counting crossings, as follows:
\begin{eqnarray*}
2A(a) &=& \diagup \hspace{-.34cm}{\color{orange2}{\nwarrow}}(a) + \diagdown \hspace{-.34cm}{\color{orange2}{\swarrow}}(a)
-\diagup \hspace{-.34cm}{\color{orange2}{\searrow}}(a) - \diagdown \hspace{-.34cm}{\color{orange2}{\nearrow}}(a),\\
M(a) & =& \diagup \hspace{-.32cm}\diagdown(a) - \diagup \hspace{-.34cm}{\color{orange2}{\searrow}}(a) -\diagdown \hspace{-.34cm}{\color{orange2}{\nearrow}}(a).
\end{eqnarray*}
Further, 
\begin{eqnarray*}
A(U_i a) & =& A(a) -1,\\
M(U_i a) &=& M(a)-2.
\end{eqnarray*}

Setting all $U_i$ to zero, we get a bigraded quotient algebra $\ah{P}= \am{P}/(U_i=0)$ over $\Ft$. 

\subsection{The bimodule for a tangle}\label{ssec:mod}

\begin{figure}[h]
\centering
    \includegraphics[scale=.82]{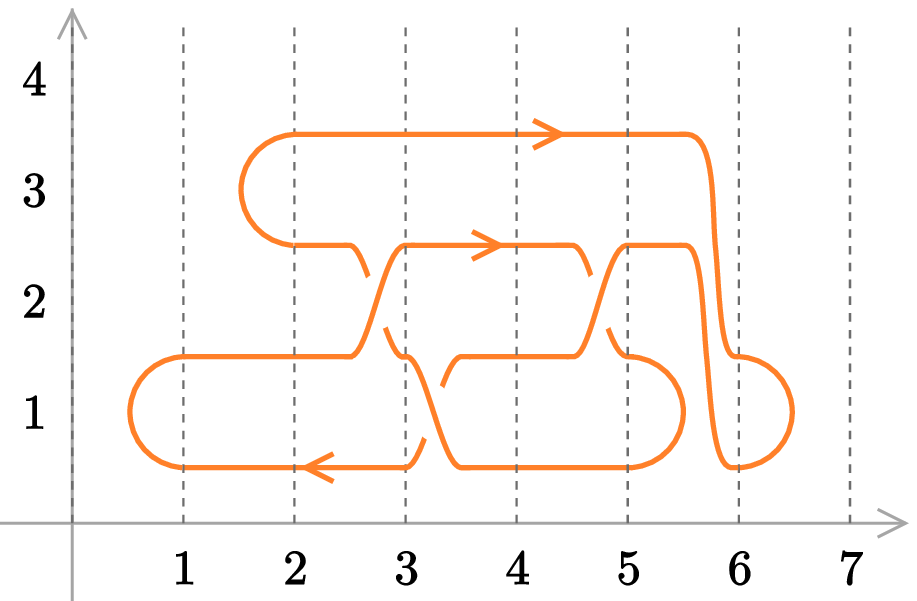} 
      \vskip .1cm
       \caption{The trefoil, decomposed as $\Tdec = (\T_1, \ldots, \T_7)$. The elementary tangles $\T_1$ and $\T_2$ are examples of cups, $\T_3$, $\T_4$, and $\T_5$ are crossings, $\T_6$ and $\T_7$ are caps. }\label{fig:tdec}
\end{figure}

Let $\Tdec$ be a decomposition of a tangle $\T$ into elementary tangles (crossings, cups, caps, or straight strands), as in Figure \ref{fig:tdec}. To $\Tdec$ we can associate a $\tDA$ structure $\ctm(\Tdec)$ over $(\am{-\bdy^L\T}, \am{\bdy^R\T})$, by looking at sequences of strand diagrams, or by looking at a plumbing of bordered grid diagrams. The two constructions are equivalent, despite using two seemingly different languages. We describe them in parallel. 

\subsubsection{The module as strand diagrams}\label{sssec:strands}

Let $\Tdec = (\T_1, \ldots, \T_k)$ be a decomposition of a tangle $\T$ into elementary tangles. Draw the projection of $\T$ in $[0,k]\times \R$ in orange, so that each $\T_i$ is in $[i-1, i]\times \R$, and:
\begin{itemize}
\item Cups and caps look like right-opening and left-opening semicircles of radius $1/2$, respectively.
 \item The two strands  at each crossing are monotone with respect to the $y$-axis; if the strand with the higher slope goes over the strand with the lower slope, then they cross in $(i-\frac{1}{2},i)\times \R$, otherwise they cross in $(i-1, i-\frac{1}{2})\times \R$. (So one can recover the type of crossing from its coordinates.)
 \item All remaining strands are monotone with respect to the $x$-axis, and don't intersect other strands. 
 \end{itemize}
See Figure \ref{fig:tdec}.

\begin{figure}[h]
\centering
  \labellist
  \pinlabel \small{$a_0^0$} at 7 15
  \pinlabel \small{$a_1^0$} at 7 47
  \pinlabel \small{$a_2^0$} at 7 79
  \pinlabel \small{$a_3^0$} at 7 111
    \pinlabel \small{$b_0^1$} at 24 15
    \pinlabel \small{$b_1^1$} at 24 47
        \pinlabel \small{$b_2^1$} at 24 79
    \pinlabel \small{$b_3^1$} at 24 111
  \pinlabel $A_0$ at 12 125
  \pinlabel $B_1$ at 30 125
  \pinlabel $A_1$ at 47 125
  \pinlabel $B_2$ at 62 125
  \pinlabel $A_2$ at 77 125
  \pinlabel $B_3$ at 92 125
  \pinlabel $A_3$ at 109 125
  \pinlabel $B_4$ at 125 125
  \pinlabel $A_4$ at 142 60
  \pinlabel $\color{gray}{1}$ at 45 4
  \pinlabel $\color{gray}{2}$ at 77 4
  \pinlabel $\color{gray}{3}$ at 109 4
  \pinlabel $\color{gray}{4}$ at 141 4
  \endlabellist
    \includegraphics[scale=1.24]{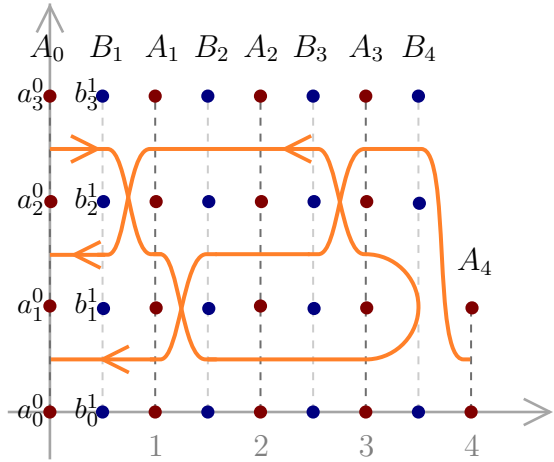} 
    \vskip .2cm
       \caption{The tangle from Figure \ref{fig:tangle} decomposed as $\Tdec = (\T_1, \T_2, \T_3, \T_4)$ (three crossings and a cap), along with the sets of points $A_0, \ldots, A_4$ (in dark red), and $B_1, \ldots, B_4$ (in dark blue). The points in $A_0$ and $B_1$ are labeled.}\label{fig:tdecAB}
\end{figure}

For $0\leq i\leq k$, pick one point in each segment of $(\{i\}\times \R)\setminus \T$, and label them $a_0^i, \ldots, a_{s_i}^i$, indexed by relative height. Let $A_i = \{a_0^i, \ldots, a_{s_i}^i\}$. Similarly, for $1\leq i\leq k$, pick one point in each segment of $(\{i-\frac{1}{2}\}\times \R)\setminus \T$, and label them $b_0^i, \ldots, b_{t_i}^i$, indexed by relative height. Let $B_i = \{b_0^i, \ldots, b_{t_i}^i\}$. See Figure \ref{fig:tdecAB}.

Let $\gen(\Tdec)$ be the set of sequences of partial permutations $x_1^-:A_0\to B_1$, $x_1^+:B_1\to A_1$, ...,  $x_k^+:B_k\to A_k$ such that $\mathrm{Image}(x_i^-) = B_i\setminus \mathrm{Domain}(x_i^+)$ and $\mathrm{Image}(x_{i-1}^+) = A_i\setminus \mathrm{Domain}(x_i^-)$. We can represent an element $\x = (x_1^-, \ldots, x_k^+)\in \gen(\Tdec)$ as a sequence of \emph{strand diagrams} --  connect each point in a domain to its image by a black strand that is monotone with respect to the $x$-axis, so that there are no triple intersections between strands of any color (black or orange). Strand diagrams are again considered up to planar isotopy (fixing the endpoints of strands) and Reidemeister III moves.

Let $\{T_1, \ldots, T_n\}$ be the set of segments of $\T\cap (\bigcup_i [i-\frac{1}{2}, i]\times \R)$ oriented left-to-right and segments of $\T\cap (\bigcup_i [i, i+\frac{1}{2}]\times \R)$ oriented right-to-left.
Let  $\ctm(\Tdec)$ be the vector space generated over $\Ft[U_1, \ldots, U_n]$ by $\gen(\Tdec)$. 
This space has an Alexander grading in $\frac{1}{2}\Z$ and a Maslov grading in $\Z$, defined on a generator $\x =(x_1^-, \ldots, x_k^+)$ as follows. Define a function $A$ on partial bijections by 
\begin{eqnarray*}
2A(x) &=&  \diagup \hspace{-.375cm}{\color{orange2}{\nwarrow}}(x) + \diagdown \hspace{-.375cm}{\color{orange2}{\swarrow}}(x)
-\diagup \hspace{-.37cm}{\color{orange2}{\searrow}}(x) - \diagdown \hspace{-.37cm}{\color{orange2}{\nearrow}}(x)
+{\color{orange2}{\searrow}} \hspace{-.38cm}{\color{orange2}{\nearrow}}(x)
-{\color{orange2}{\nwarrow}} \hspace{-.38cm}{\color{orange2}{\swarrow}}(x)
-{\color{orange2}{\leftarrow}}(x),
\end{eqnarray*}
for $x\in \{x_1^-, \ldots, x_k^+\}$, and define the Alexander grading of $\x$ by  $A(\x) = A(x_1^-)+ \cdots + A(x_k^+)$. 

Define a function $M$ on  partial bijections by
\begin{eqnarray*}
M(x_i^-)  &=& -\diagup \hspace{-.34cm}\diagdown(x_i^-) + \diagup \hspace{-.375cm}{\color{orange2}{\nwarrow}}(x_i^-) + \diagdown \hspace{-.375cm}{\color{orange2}{\swarrow}}(x_i^-)
-{\color{orange2}{\swarrow}} \hspace{-.38cm}{\color{orange2}{\nwarrow}}(x_i^-)
-{\color{orange2}{\leftarrow}}(x_i^-),\\
M(x_i^+) &=& \diagup \hspace{-.34cm}\diagdown(x_i^+) - \diagup \hspace{-.37cm}{\color{orange2}{\searrow}}(x_i^+) -\diagdown \hspace{-.37cm}{\color{orange2}{\nearrow}}(x_i^+)
+{\color{orange2}{\searrow}} \hspace{-.38cm}{\color{orange2}{\nearrow}}(x_i^+),
\end{eqnarray*}
and define the Maslov grading of $\x$ by $M(\x) = M(x_1^-) + \cdots + M(x_k^+)$.

Further, 
\begin{eqnarray*}
A(U_i \x) & =& A(\x) -1,\\
M(U_i \x) &=& M(\x)-2.
\end{eqnarray*}

\begin{figure}[h]
\centering
  \includegraphics[scale=.75]{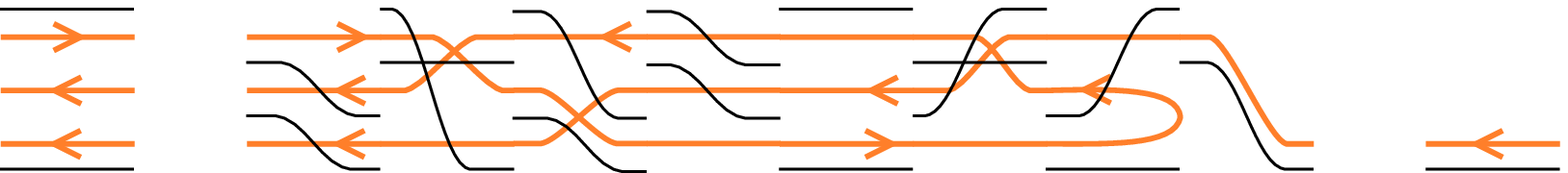} 
    \vskip .2cm
       \caption{Middle: a generator $\x \in \gen(\Tdec)$ for the decomposition $\x$ in Figure \ref{fig:tdecAB}. Left: the idempotent $\eeD{L}(\x)$. Right: the idempotent $\eeA{R}(\x)$.}\label{fig:gen_strands}
\end{figure}

\begin{example*} \emph{The generator in Figure \ref{fig:gen_strands} has Alexander grading $-1$ and Maslov grading  $-1$.}
\end{example*}

Think of $\am{-\bdy^L\T}$ and $\am{\bdy^R\T}$ as generated by partial permutations on  $A_0$ and  $A_k$ instead of on  $[|-\bdy^L\T|]$ and $[|\bdy^R\T|]$ (the goal is to soon define the $\tDA$ structure maps graphically, via strand diagrams). 
 We give $\ctm(\Tdec)$ the structure of a left-right bimodule over $(\im{-\bdy^L\T}, \im{\bdy^R\T})$ by defining
\begin{equation*}
\eee{\sss} \x \eee{\ttt} = \begin{cases}
\x &  \sss = A_0\setminus \mathrm{Domain}(\x), \ttt = A_k\cap \mathrm{Image}(\x) \\
0 & \text{otherwise}.
\end{cases} 
\end{equation*}
In other words, $\eee{\sss}$ acts on the left by identity when $\sss$ is the elements of $A_0$ unoccupied by $\x$ (denote this idempotent by $\eeD{L}(\x)$), and by zero otherwise, and  $\eee{\ttt}$ acts on the right by identity when $\ttt$ is the elements of $A_k$ occupied by $\x$ (denote this idempotent by $\eeA{R}(\x)$), and by zero otherwise.
See Figure \ref{fig:gen_strands}.

We next define three maps $\bdy_+, \bdy_-, \bdy_m: \ctm(\Tdec)\to \ctm(\Tdec)$.

\begin{figure}[h]
\centering
  \labellist
  \pinlabel $\color{orange2}{T_i}$ at -8 43
  \pinlabel $=U_i$ at 61 49
  \pinlabel $\color{orange2}{T_i}$ at -8 14
  \pinlabel $=U_i$ at 61 10
  \pinlabel $=0$ at 227 50
  \pinlabel $=0$ at 227 11
  \pinlabel $=0$ at 396 50
  \endlabellist
\hspace{0.7em}   \includegraphics[scale=.8]{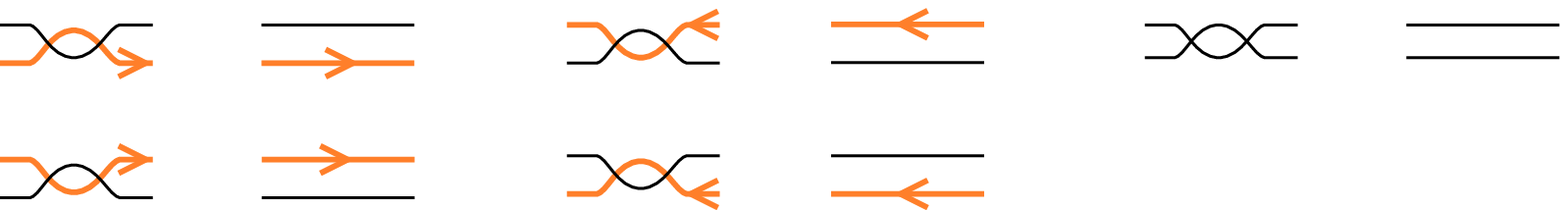} 
  \vskip .2cm
       \caption{Relations for $\bdy_+$. The relations in the second and third column just say that those moves are not allowed.}\label{fig:Reid}
\end{figure}

For a generator $\x$, $\bdy_+(\x)$ is the sum of all elements of $\ctm(\Tdec)$ obtained by \emph{smoothing a black-black  crossing of $\x$ contained in $(\bigcup_i (i-\frac{1}{2}, i)\times \R)$} and then performing necessary Reidemeister II moves to obtain a valid strand diagram, where the relations in Figure \ref{fig:Reid} are applied to the Reidemeister II moves. Graphically, $\bdy_+$ is the same as the differential on the algebra. Extend $\bdy_+$ linearly to all of $\ctm(\Tdec)$.

\begin{figure}[h]
\centering
  \labellist
  \pinlabel $\color{orange2}{T_i}$ at -8 43
  \pinlabel $=U_i$ at 61 49
  \pinlabel $\color{orange2}{T_i}$ at -8 14
  \pinlabel $=U_i$ at 61 10
  \pinlabel $=0$ at 227 50
  \pinlabel $=0$ at 227 11
  \pinlabel $=0$ at 396 50
  \endlabellist
\hspace{0.7em}    \includegraphics[scale=.8]{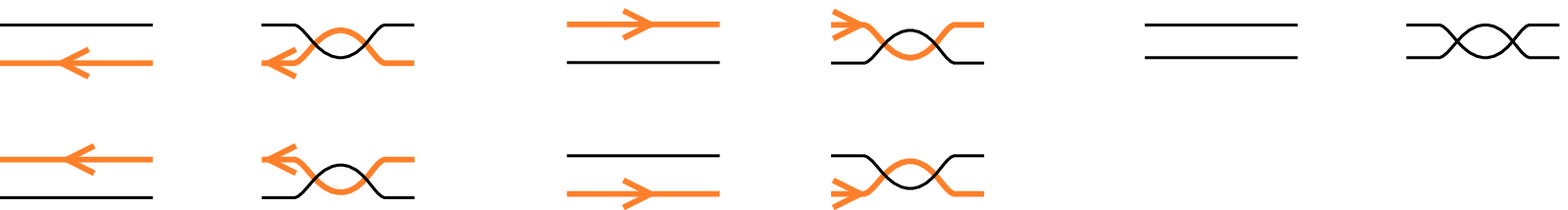} 
  \vskip .2cm
       \caption{Relations for $\bdy_-$. The relations in the second and third column just say that those moves are not allowed.}\label{fig:Reid_}
\end{figure}

For a generator $\x$,  $\bdy_-(\x)$ is the sum of all elements of $\ctm(\Tdec)$ obtained by \emph{introducing a black-black crossing to $\x$ in $(\bigcup_i (i, i+\frac{1}{2})\times \R)$}, by  performing  Reidemeister II moves if necessary to bring a pair of non-crossing black strands close together, where the  relations in Figure \ref{fig:Reid_} are applied to the Reidemeister II moves.
Extend $\bdy_-$ linearly to all of $\ctm(\Tdec)$.

\begin{figure}[h]
\centering
    \includegraphics[scale=.73]{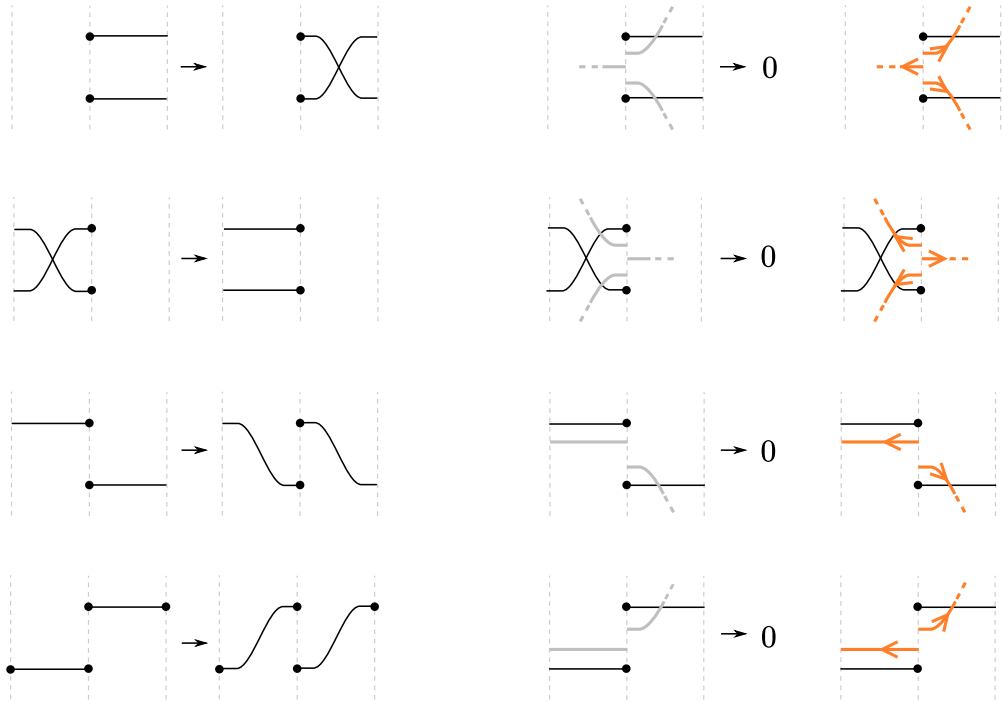} 
      \vskip .2cm
       \caption{The map $\bdy_m$ along a vertical line $L$ containing a set $B_j$. Left: The map $\bdy_m$ counts four types of exchanges of pairs of black strands with ends at $B_j$, modulo the relations depicted in the middle and right column. Middle: For a fixed pair of black strands, if there is  an orange segment oriented \emph{into} $L$ or a black strand in the same relative position to the pair as one of the displayed grey strands, the exchange is not allowed. Right: Each orange segment $T_i$ oriented \emph{away} from $L$ and in the same relative position to the pair as one of the displayed grey strands in the middle column results in a factor of $U_i$ for the resulting exchange. The map $\bdy_m$ along a line containing a set $A_j$ can be described graphically  by reflecting each digram in this figure along a vertical line, and switching orientations on the orange segments in the relations.} \label{fig:d_m_minus}
\end{figure}

For a generator $\x$,  $\bdy_m(\x)$ is the sum of all elements of $\ctm(\Tdec)$ obtained by picking a pair of points in a given $A_j$ or a given $B_j$,  and \emph{exchanging these two ends of the corresponding pair of black strands of $\x$}, subject to the relations in Figure \ref{fig:d_m_minus}. 
Extend $\bdy_m$ linearly to all of $\ctm(\Tdec)$.

\begin{figure}[h]
\centering
  \labellist
  \pinlabel $\color{orange2}{T_3}$ at -8 123
  \pinlabel $\color{orange2}{T_4}$ at -8 140
  \pinlabel $\color{orange2}{T_1}$ at 88 90
  \pinlabel $\color{orange2}{T_2}$ at 88 125
  \pinlabel $U_3$ at 125 118
   \pinlabel $+U_3$ at 234 118
   \pinlabel $+U_3$ at 234 33
   \pinlabel $+U_2$ at 348 118
   \pinlabel $+$ at 348 33
  \pinlabel $+U_2U_3$ at 113 33
  \endlabellist
\hspace{0.7em}    \includegraphics[scale=.8]{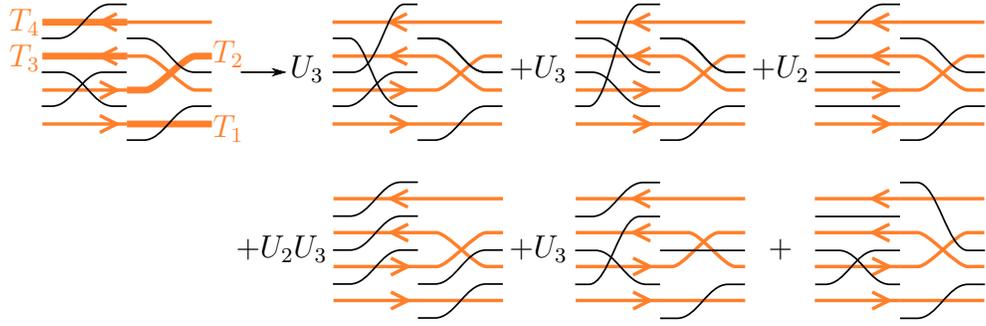} 
  \vskip .2cm
       \caption{An example of the map $\bdy_+ + \bdy_-+\bdy_m$, applied to a generator $\x$ for a tangle decomposition consisting of a single crossing. The segments $T_i$ are drawn thicker in the starting generator. The first two terms on the right hand side correspond to $\bdy_-$, the last four to $\bdy_m$, and $\bdy_+$ is zero.}\label{fig:d_minus}
\end{figure}

See Figure \ref{fig:d_minus} for an example of the map $\bdy_+ + \bdy_-+\bdy_m$.

\begin{figure}[h]
\centering
  \labellist
  \pinlabel $\delta_1^1(\x)=$ at 42 355
  \pinlabel $+$ at 60 277
  \pinlabel $\otimes$ at 130 355
  \pinlabel $\otimes$ at 130 277
  \pinlabel $+$ at 60 200
  \pinlabel $\otimes$ at 130 200
  \pinlabel $+$ at 60 124
  \pinlabel $\otimes$ at 130 124
  \pinlabel $\delta_2^1\Big(\x\otimes\hspace{33pt}\Big)=$ at 9 24
  \pinlabel $\otimes$ at 130 24
  \endlabellist
\hspace{3em}    \includegraphics[scale=.75]{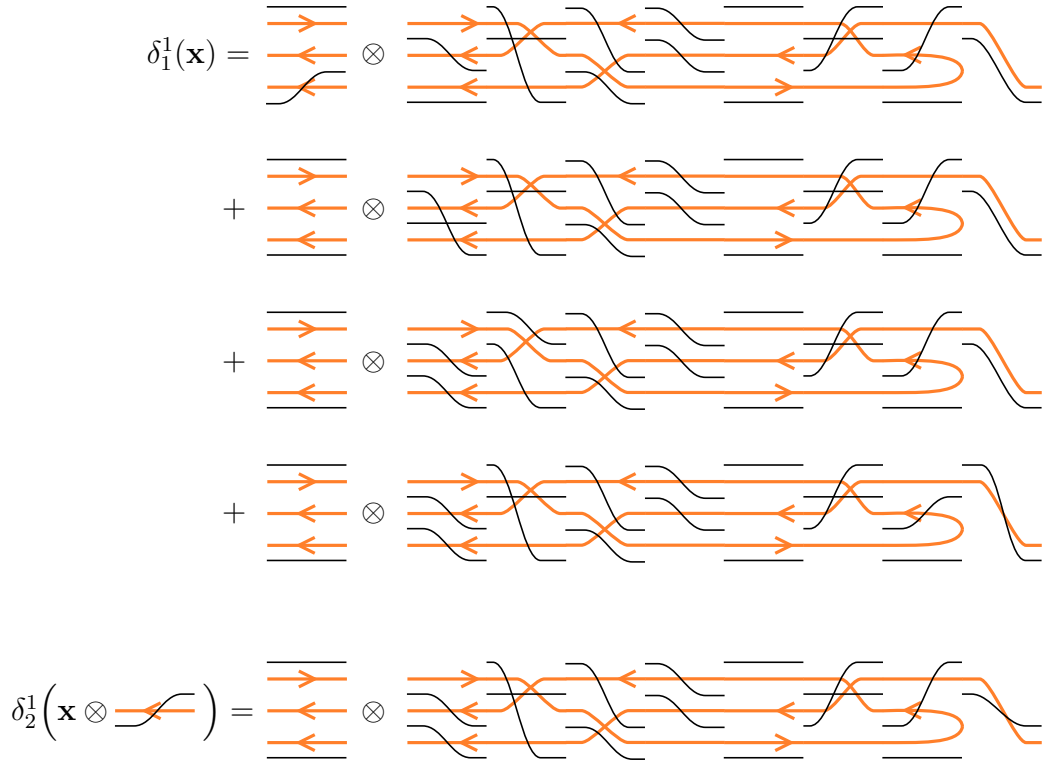} 
  \vskip .2cm
       \caption{The $\tDA$ structure maps for the generator $\x$ in Figure \ref{fig:gen_strands}.}\label{fig:DA_strands}
\end{figure}

Next, we define a map $m_2: \ctm(\Tdec)\otimes_{\im{\bdy^R\T}}\am{\bdy^R\T}\to \ctm(\Tdec)$. Here,  each variable $U_i$ in $\ctm(\Tdec)$ corresponding to a segment $T_i$ of the tangle in the rightmost piece $[k-\frac{1}{2},k]\times \R$ oriented left-to-right is identified with the variable $U_j$ in $\am{\bdy^R\T}$ corresponding to $(T_i\cap (\{k\}\times \R))\subset \bdy^R\T$.
On generators, define  $m_2((x_1^-, \ldots, x_k^-,  x_k^+), a) = (x_1^-, \ldots, x_k^-, x_k^+a)$ where $x_k^+a$ is given by concatenating $x_k^+$ to the left with $a$ to the right, modulo the relations in Figure \ref{fig:Reid}. If $x_k^+$ and $a$ cannot be concatenated, then  $m_2((x_1^-, \ldots, x_k^-,  x_k^+), a) = 0$. See the bottom of Figure \ref{fig:DA_strands}. Extend $m_2$ linearly to all of $\ctm(\Tdec)\otimes_{\im{\bdy^R\T}}\am{\bdy^R\T}$.

Last, define a map $\delta^L:  \ctm(\Tdec)\to \am{-\bdy^L\T}\otimes_{\im{-\bdy^L\T}} \ctm(\Tdec)$. For a generator $\x$, $\delta^L(\x)$ is given by gluing a diagram for  $\eeD{L}(\x)$ to the left of a diagram for $\x$, and them applying the same exchange map as $\bdy_m$ to the gluing line. See the top of Figure \ref{fig:DA_strands}.

The above maps can be combined to define a $\tDA$ structure. 

\begin{definition}
We give the $(\im{-\bdy^L\T},\im{\bdy^R\T})$ bimodule  $\ctm(\Tdec)$ the structure of a type $\tDA$ bimodule over $(\am{-\bdy^L\T},\am{\bdy^R\T})$ using the following structure maps: Define
$$\delta_1^1: \ctm(\Tdec)\to  \am{-\bdy^L\T}\otimes_{\im{-\bdy^L\T}}\ctm(\Tdec)$$
on generators by 
$$\delta_1^1(\x) = \eeD{L}(\x)\otimes (\bdy_+ + \bdy_- + \bdy_m)(\x) + \delta^L(\x), $$
 define 
$$\delta_2^1: \ctm(\Tdec)\otimes_{\im{\bdy^R\T}}\am{\bdy^R\T}\to  \am{-\bdy^L\T}\otimes_{\im{-\bdy^L\T}}\ctm(\Tdec)$$
on generators by 
$$\delta_2^1(\x \otimes a) = \eeD{L}(\x)\otimes m_2(\x, a),$$
and define $\delta_i^1=0$ for $i>2$.
\end{definition}

See \cite{pv} for a proof that this is indeed a type $\tDA$ structure, $\delta_1^1$ lowers the Maslov grading by one, and preserves the Alexander grading. 

For example, Figure \ref{fig:DA_strands} shows the type $\tDA$ structure maps applied to the generator from Figure \ref{fig:gen_strands}.

Note that $\ctm(\Tdec)$ splits as the direct sum $\bigoplus\ctm_i(\Tdec)$, where $\ctm_i(\Tdec)$ is the structure generated by elements with $i$ black strands in the rightmost piece. 

In \cite{pv}, we prove the first part of Theorem \ref{thm:main} by observing that concatenating two tangle decompositions  corresponds to taking a box tensor product of their respective $\tDA$ structures.

\subsubsection{The module as bordered grid diagrams}\label{sssec:grids}

\begin{figure}[h]
\centering
  \labellist
  \pinlabel \small{$\color{red}{\alpha_0^0}$} at 22 23
   \pinlabel $\rotatebox{90}{\scalebox{1.5}[1.0]{\color{red}{\ldots}}}$ at 19 45
  \pinlabel \small{$\color{red}{\alpha_3^0}$} at 20 65
  \pinlabel \small{$\color{blue}{\beta_0^1}$} at 64 4
  \pinlabel $\rotatebox{90}{\scalebox{1.4}[1.0]{\color{blue}{\ldots}}}$ at 64 20
  \pinlabel \small{$\color{blue}{\beta_3^1}$} at 64 35
  \pinlabel \small{$\color{red}{\alpha_0^1}$} at 110 24
  \pinlabel $\rotatebox{90}{\scalebox{1.5}[1.0]{\color{red}{\ldots}}}$ at 110 44
  \pinlabel \small{$\color{red}{\alpha_3^1}$} at 110 64
  \pinlabel \small{$\color{blue}{\beta_0^2}$} at 155 4
  \pinlabel $\rotatebox{90}{\scalebox{1.4}[1.0]{\color{blue}{\ldots}}}$ at 155 20
  \pinlabel \small{$\color{blue}{\beta_3^2}$} at 155 35
    \pinlabel \small{$\color{red}{\alpha_0^2}$} at 201 24
  \pinlabel $\rotatebox{90}{\scalebox{1.5}[1.0]{\color{red}{\ldots}}}$ at 201 44
  \pinlabel \small{$\color{red}{\alpha_3^2}$} at 201 64
    \pinlabel \small{$\color{blue}{\beta_0^3}$} at 246 4
   \pinlabel $\rotatebox{90}{\scalebox{1.4}[1.0]{\color{blue}{\ldots}}}$ at 246 20
  \pinlabel \small{$\color{blue}{\beta_3^3}$} at 246 35
      \pinlabel \small{$\color{red}{\alpha_0^3}$} at 292 24
  \pinlabel $\rotatebox{90}{\scalebox{1.5}[1.0]{\color{red}{\ldots}}}$ at 292 44
  \pinlabel \small{$\color{red}{\alpha_3^3}$} at 292 64
    \pinlabel \small{$\color{blue}{\beta_0^4}$} at 331 11
    \pinlabel \small{$\color{blue}{\beta_1^4}$} at 331 23
  \pinlabel \small{$\color{blue}{\beta_2^4}$} at 331 35
      \pinlabel \small{$\color{red}{\alpha_0^4}$} at 366 34
  \pinlabel \small{$\color{red}{\alpha_1^4}$} at 366 55
  \endlabellist
    \includegraphics[scale=1.04]{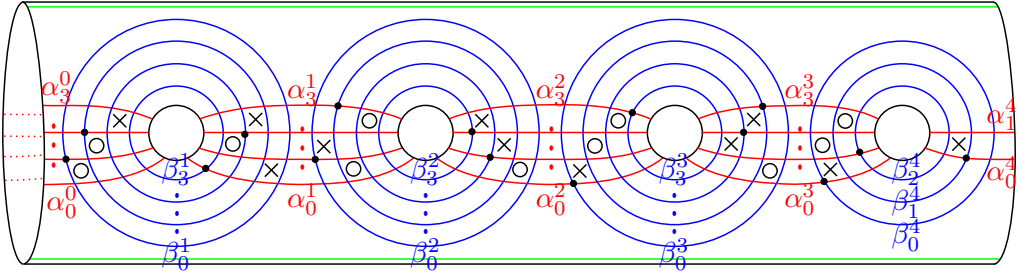} 
      \vskip .2cm
       \caption{The Heegaard diagram and generator corresponding to Figure \ref{fig:gen_strands}.}\label{fig:gen_hd}
\end{figure}

To a tangle decomposition $\Tdec = (\T_1, \ldots, \T_k)$ one can also associate a  bordered Heegaard diagram $\mathcal H(\Tdec)$, as follows. Start with a genus $k$ surface $\Sigma$ with two boundary components $\bdy^L\Sigma$ and $\bdy^R\Sigma$, and draw parallel $\beta$ circles, one circle $\beta_t^j$ for each $b_t^j\in B_j$, for $1\leq j\leq k$, and parallel $\alpha$ circles, one circle $\alpha_t^j$ for each $a_t^j\in A_j$, for $1\leq j\leq k-1$, as well as parallel $\alpha$ arcs, one $\alpha_t^0$ for each $a_t^0\in A_0$ with ends on $\bdy^L\Sigma$, and parallel $\alpha$ arcs, one $\alpha_t^k$ for each $a_t^k\in A_k$ with ends on $\bdy^R\Sigma$, as in Figure \ref{fig:gen_hd}. 

If there is a segment of the tangle oriented left-to-right, respectively right-to-left,  running from somewhere between  $b_s^j$ and $b_{s+1}^j$ to somewhere between  $a_t^j$ and $a_{t+1}^j$, place an $O$, respectively $X$, on the Heegaard diagram, so that it is contained  in the annulus bounded by $\beta_s^j$ and $\beta_{s+1}^j$, as well as in  the annulus bounded by $\alpha_t^j$ and $\alpha_{t+1}^j$.

If there is a segment of the tangle oriented  left-to-right, respectively right-to-left, running from somewhere between  $a_t^j$ and $a_{t+1}^j$ to somewhere between $b_s^{j+1}$ and $b_{s+1}^{j+1}$, place an $X$, respectively $O$, on the Heegaard diagram, so that it is contained  in the annulus bounded by $\alpha_t^j$ and $\alpha_{t+1}^j$, as well as in  the annulus bounded by $\beta_s^{j+1}$ and $\beta_{s+1}^{j+1}$.

One can see the tangle by connecting $X$s to $O$s by arcs away from the $\beta$ curves and pushing those arcs slightly above the Heegaard surface, and connecting $O$s to $X$s by arcs away from the $\alpha$ curves, as well as $O$s to points on $\bdy\Sigma$, and points on $\bdy\Sigma$ to $X$s away from the $\alpha$ curves.

Define the \emph{generators} of $\mathcal H(\Tdec)$ to be sets of intersection points of $\alpha$ and $\beta$ curves, so that there is exactly one point on each $\beta$ circle and on each $\alpha$ circle, and at most one point on each $\alpha$ arc. Note that these are in one-to-one correspondence with generators in $\gen(\Tdec)$ (a strand connecting $b_s^i$ to  $a_t^j$, if those two points are in adjacent sets, corresponds to the point $\beta_s^i\cap \alpha_t^j$). Grade the generators  of $\mathcal H(\Tdec)$ the same as their corresponding generators in $\gen(\Tdec)$.

\begin{figure}[h]
\centering
    \includegraphics[scale=.92]{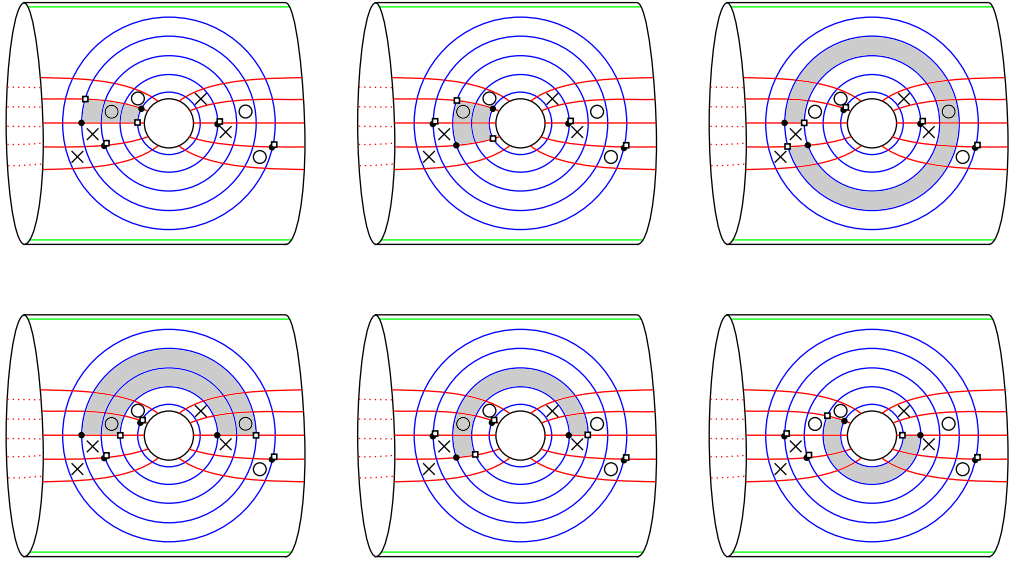} 
      \vskip .2cm
       \caption{The six terms in $\bdy(\x)$, for the generator from Figure \ref{fig:d_minus}, along with the respective rectangles, in order consistent with Figure \ref{fig:d_minus}.}\label{fig:d_minus_hd}
\end{figure}

The map $(\bdy_+ + \bdy_-+\bdy_m)(\x)$ corresponds to the map $\bdy$ that counts internal rectangles in $\mathcal H(\Tdec)$ that are empty (the interior contains no intersection points of $\x$ and no $X$s), so that crossing an $O$ corresponding to a segment $T_i$ in the tangle decomposition results in multiplication by $U_i$. See Figure \ref{fig:d_minus_hd} for an example.

\begin{figure}[h]
\centering
    \includegraphics[scale=.78]{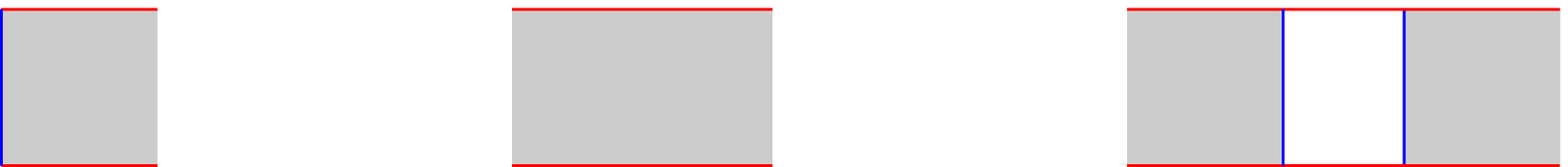} 
      \vskip .2cm
       \caption{The three types of rectangles for $\bdy^L$.}\label{fig:D-rect}
\end{figure}

The map $\delta^L$ corresponds to a map $\bdy^L$ that counts the following types of embedded rectangles that intersect $\bdy^L\Sigma$:

\begin{enumerate}
\item A rectangle $r$ whose oriented boundary follows an arc of $\alpha_i^0$, then an arc of $\beta_m^1$, then an arc of $\alpha_j^0$, then an arc of $\bdy^L\Sigma$. Given generators $\x$ and $\y$, $r$ \emph{connects $\x$ to  $\y$} if $\alpha_j^0\cap \beta_m^1 = \x\setminus \y$ and $\alpha_i^0\cap \beta_m^1 = \y\setminus \x$. Define $a^{r, \x, \y}\in \am{-\bdy^L\T}$ as the bijection with domain $\{t|\alpha_t^0\cap \x = \emptyset\}$ that sends $i$ to $j$, and is the identity everywhere else. Define $U^r$ as the product of all variables $U_s$ with corresponding $O_s$ in the interior of $r$.
\item A rectangle $r$ whose boundary consists of two complete arcs $\alpha_i^0$ and $\alpha_j^0$ and two arcs in $\bdy^L\Sigma$.  Given generators $\x$ and $\y$, $r$ \emph{connects $\x$ to  $\y$} if $\x = \y$, and $\alpha_i^0$ and $\alpha_j^0$ are not occupied by $\x=\y$.  Define $a^{r, \x, \y}\in \am{-\bdy^L\T}$ as the bijection with domain $\{t|\alpha_t^0\cap \x = \emptyset\}$ that sends $i$ to $j$, $j$ to $i$, and is the identity everywhere else. Define $U^r$ as the product of all variables $U_s$ with corresponding $O_s$ in the interior of $r$.
\item For $i<j$ and $m<n$, the union $r$ of two disjoint rectangles of the first type, such that one has boundary on  $\alpha_i^0$, $\beta_m^1$,  $\alpha_j^0$, $\bdy^L\Sigma$, and the other has boundary on $\alpha_j^0$, $\beta_n^1$,  $\alpha_i^0$, $\bdy^L\Sigma$.  Given generators $\x$ and $\y$, $r$ \emph{connects $\x$ to  $\y$} if $\{\alpha_j^0\cap \beta_m^0, \alpha_i^0\cap \beta_n^1\} = \x\setminus \y$ and $\{\alpha_i^0\cap \beta_m^1, \alpha_j^0\cap \beta_n^1\}  = \y\setminus \x$. Define $a^{r, \x, \y}\in \am{-\bdy^L\T}$ as the identity bijection with domain $\{t|\alpha_t^0\cap \x = \emptyset\}$. Define $U^r$ as the product of all variables $U_s$ with corresponding $O_s$ in the interior of $r$, and all variables $U_t$ corresponding to $+$ points in $-\bdy^L\T$ that are above the $i^{\mathrm{th}}$ and below the $j^{\mathrm{th}}$ point.
\end{enumerate}
See Figure \ref{fig:D-rect}. Note that a given rectangle may connect more than one pair of generators. 

\begin{figure}[h]
\centering
    \includegraphics[scale=.92]{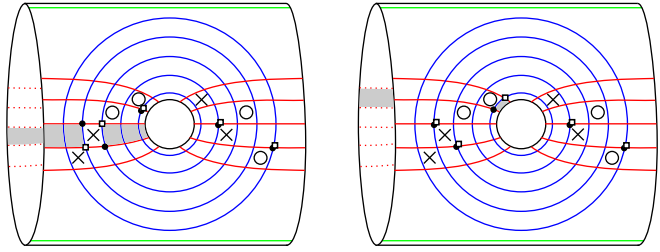} 
      \vskip .2cm
       \caption{The two terms in  $\bdy^L(\x)$, for the generator from Figure \ref{fig:d_minus} and the respective rectangles, in order consistent with Figure \ref{fig:d_minus}.}\label{fig:left_hd}
\end{figure}

 The first two types of rectangles are  \emph{empty} if their interior contains no intersection points of $\x$ and no $X$s. The third type is \emph{empty} if, in addition to its interior containing no intersection points of $\x$ and no $X$s, the interior of the internal rectangle with boundary on $\alpha_i^0$, $\beta_n^1$,  $\alpha_j^0$, $\beta_m^1$ contains $j-i-1$ $X$s and $j-i-1$ points of $\x$.
$$\bdy^L(\x) = \sum_{\y - \textrm{generator}}\sum_{r - \textrm{empty rectangle from $\x$ to $\y$}}U^ra^{r, \x, \y}\y$$
See Figure \ref{fig:left_hd} for an example of $\bdy^L$. 

The map $m_2$ corresponds to counting sets of embedded rectangles  $r = \{r_1, \ldots, r_l\}$ that intersect $\bdy^R\Sigma$, i.e. rectangles $r_s$ whose oriented boundary follows an arc of $\alpha_{i_s}^k$, then an arc of $\beta_{m_s}^k$, then an arc of $\alpha_{j_s}^k$, then an arc of $\bdy^R\Sigma$. 
 The set $r$ \emph{connects a generator $\x$ to a generator $\y$} if $\{\alpha_{j_s}^k\cap \beta_{m_s}^k| 1\leq s\leq l\} = \x\setminus \y$ and $\{\alpha_{i_s}^k\cap \beta_{m_s}^k| 1\leq s\leq l\} = \y\setminus \x$. Define $a^{r, \x, \y}\in \am{\bdy^R\T}$ as the bijection with domain $\{i|\alpha_i^k\cap \x \neq \emptyset\}$ that sends $j_s$ to $i_s$ for $1\leq s\leq l$, and is the identity everywhere else. Again note that a given set of rectangles may connect more than one pair of generators. Define $U^{r_s}$ as the product of all variables $U_t$ with corresponding $O_t$ in the interior of $r_s$.

\begin{figure}[h]
\centering
    \includegraphics[scale=.72]{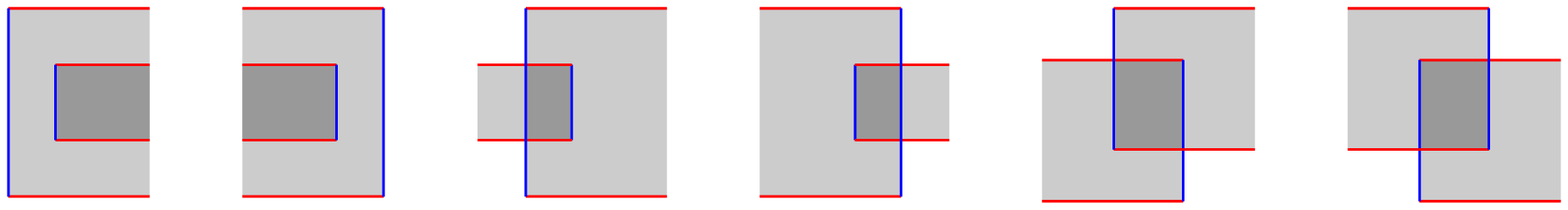} 
      \vskip .2cm
       \caption{Forbidden pairs of rectangles for $\bdy^R$.}\label{fig:forbidden_rect}
\end{figure}

A set $r = \{r_1, \ldots, r_n\}$ connecting $\x$ to $\y$ is \emph{allowed} if there are no $X$s and no points in $\x\cap \y$ in the interior of $r_s$, for $1\leq s\leq l$, and no two rectangles are in relative position as in Figure \ref{fig:forbidden_rect}. 

\begin{figure}[h]
\centering
    \includegraphics[scale=.92]{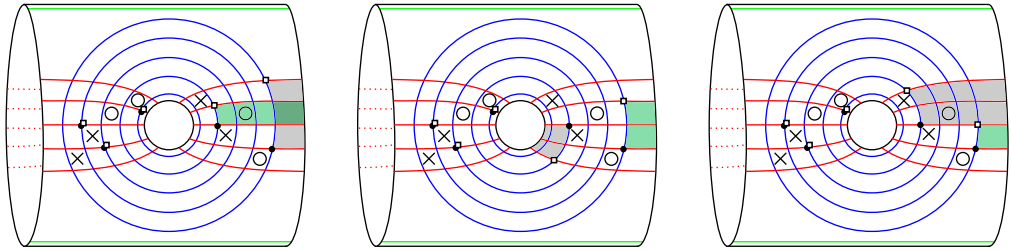} 
      \vskip .2cm
       \caption{Some sets rectangles contributing to $m_2(\x,a)$, for the generator $\x$ from Figure \ref{fig:d_minus}, and various choices of $a$.}\label{fig:right_hd}
\end{figure}

Note that for a fixed generator $\x$ and algebra generator $a$, there is at most one $\y$ and at most one $r$ as above. Thus, if there is no generator $\y$ and no allowed set of rectangles $r$ from $\x$ to $\y$ with $a^{r, \x, \y}=a$, we define $m_2(\x, a)=0$. If there are such $\y$ and $r=\{r_1, \ldots, r_l\}$, we define 
$$m_2(\x,a) =\prod_{1\leq s\leq l}U^{r_s} \y.$$
See Figure \ref{fig:right_hd} for an example. 

\begin{figure}[h]
\centering
  \labellist
  \pinlabel $1$ at 20 34
  \pinlabel $2$ at 27 44
  \pinlabel $3$ at 83 39
  \pinlabel $4$ at 330 27
  \pinlabel $5$ at 363 44
  \endlabellist
    \includegraphics[scale=1.07]{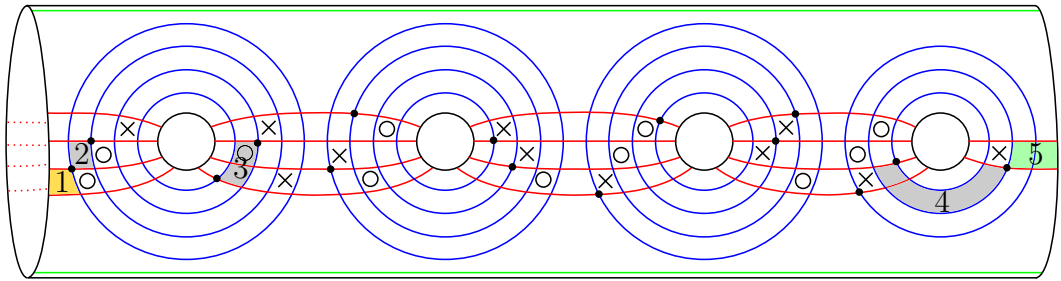} 
      \vskip .2cm
       \caption{The $\tDA$ structure maps applied to the generator from Figure \ref{fig:gen_strands} count rectangles. The indices of the five shaded rectangles correspond to the rows in Figure \ref{fig:DA_strands}.}\label{fig:DA_hd}
\end{figure}

Figure \ref{fig:DA_hd} shows the rectangles that contribute to the structure maps for the generator from Figure \ref{fig:gen_strands}.

The first part of Theorem \ref{thm:main} is proven in \cite{pv} alternatively by observing that gluing the Heegaard diagrams for the two decompositions results in a Heegaard diagram for the concatenation of the decompositions, and algebraically corresponds to taking a box tensor product as well. 

\subsection{One-sided modules and chain complexes}

When $\bdy^L\T = \emptyset$,  the left algebra is just $\am{-\bdy^L\T}\cong \Ft\oplus \Ft$, and we can think of $\ctm(\Tdec)$ as a right type $A$ structure. Similarly, when $\bdy^R \T = \emptyset$, we have $\am{\bdy^L\T}\cong \Ft\oplus \Ft$, and we can think of $\ctm(\Tdec)$ as a left type $D$ structure. When $\T$ is a closed link, $\ctm(\Tdec)$ is just a chain complex. 
In any of these cases, in \cite{pv} we make some non-canonical choices to define one-sided bordered Heegaard diagrams for tangles, and closed Heegaard diagrams for tangles, as follows.

When $\bdy^L\T = \emptyset$, the only nonzero summands of $\ctm(\Tdec)=\bigoplus\ctm_i(\Tdec)$ are $\ctm_{\frac{1}{2}|\bdy^R \T|}(\Tdec)$ and $\ctm_{\frac{1}{2}|\bdy^R \T|+1}(\Tdec)$. 
We can modify the Heegaard diagram for $\Tdec$ to only have right boundary, by gluing the ``front" and ``back" of $\bdy^L\Sigma$, so that $\alpha_0^0$ becomes a circle. The combinatorics (generators and rectangle counts) of the resulting diagram are the same as for $\mathcal H(\Tdec)$ when $\alpha_0^0$ is occupied, and we get a  type $A$ structure  that is exactly $\ctm_{|\bdy^R \T|}(\Tdec)$ (we can think of $\ctm_{|\bdy^R \T|}(\Tdec)$ as a type $A$ instead of a type $\tDA$ structure, since the left algebra for this summand is just $\Ft$).  This is the type $1$ Heegaard diagram for a tangle with only right boundary that we use in \cite{pv}. If instead we delete the new $\alpha$ circle, we obtain a diagram with corresponding type $A$ structure $\ctm_{|\bdy^R \T|+1}(\Tdec)$.

 Similarly, when $\bdy^R \T = \emptyset$ we can modify the Heegaard diagram for $\Tdec$ to only have right boundary, by gluing the ``front" and ``back" of $\bdy^R\Sigma$, so that $\alpha_0^k$ becomes a circle, and then deleting that circle. Call the resulting diagram $\mathcal H^D(\Tdec)$. The type $D$ structure corresponding to $\mathcal H^D(\Tdec)$ is $\ctm_0(\Tdec)$. This is the type $2$ Heegaard diagram for a tangle with only left boundary that we use in \cite{pv}.  If instead we leave the circle in, we get the type $D$ structure $\ctm_1(\Tdec)$.

 When $\T$ is a closed link, we can modify the diagram in both of the above ways (leaving in the closure of $\alpha_0^0$ and deleting the closure of $\alpha_0^k$) and we can place an $X$ and an $O$ in the single non-combinatorial region, since by definition of $\ctm$ that region is not considered. We get a Heegaard diagram for $\T$ union an unknot, so $\ctm_0(\Tdec)\simeq \HFKm(\T)\otimes (\Ft\oplus\Ft)$. Similarly we can delete the closure of $\alpha_0^0$ and leave the closure of $\alpha_0^k$, to see that $\ctm_1(\Tdec)\simeq \HFKm(\T)\otimes (\Ft\oplus\Ft)$ as well.
 The second part of Theorem \ref{thm:main} follows (in \cite{pv} we make the above choices and refer to what here is the summand $\ctm_0(\mathcal L)$ as all of $\ctm(\mathcal L)$).


\section{How not to compute $\HFKm$ of the unknot} 
\label{sec:computations}

We conclude this paper by working out a very small example. 
We compute the well known knot Floer homology of the unknot $\mathcal U$, by decomposing it as $\mathcal U = \T_1\circ \T_2$, where  $\T_1$ is a single cup with $\bdy^R\T_1 = (-,+)$, $\T_2$ is a single cap with $-\bdy^L\T_2 = (-,+)$, 
 computing the type $A$ structure $\T_1$ , the type $D$ structure for $\T_2$, and taking their tensor product over $\am{-,+}$.
 
\begin{figure}[h]
\centering
  \labellist
  \pinlabel \textcolor{gray}{$A$} at 30 20
  \pinlabel \textcolor{gray}{$B$} at 30 64
  \pinlabel \textcolor{gray}{$C$} at 77 29
  \pinlabel \textcolor{gray}{$D$} at 78 42
  \pinlabel \textcolor{gray}{$E$} at 54 33
  \pinlabel \textcolor{gray}{$F$} at 56 42
  \endlabellist
    \includegraphics[scale=1.25]{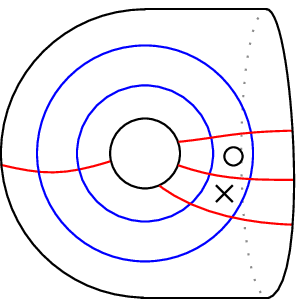} 
      \vskip .2cm
       \caption{The Heegard diagram for the single cup $\T_1$ with $\bdy^R\T_1 = (-,+)$.}\label{fig:cap_hd}
\end{figure}

 The Heegaard diagram for $\T_1$ is displayed in Figure \ref{fig:cap_hd}.
Let $a_{ij}\in \am{-,+}$ be  the algebra element consisting of a single strand connecting $i$ on the left to $j$ on the right,  and let $U_1$ be the variable corresponding to the single $O$ in the Heegaard diagram, and to the $+$ in $\bdy^R\T_1$. 

The Heegaard diagram has six generators: $x_0 = \{\alpha_0^0\cap \beta_0^1, \alpha_0^1\cap \beta_1^1\}$, $x_1 = \{\alpha_0^0\cap \beta_0^1, \alpha_1^1\cap \beta_1^1\}$, $x_2 = \{\alpha_0^0\cap \beta_0^1, \alpha_2^1\cap \beta_1^1\}$, $y_0 = \{\alpha_0^0\cap \beta_1^1, \alpha_0^1\cap \beta_0^1\}$, $y_1 = \{\alpha_0^0\cap \beta_1^1, \alpha_1^1\cap \beta_0^1\}$, $y_2 = \{\alpha_0^0\cap \beta_1^1, \alpha_2^1\cap \beta_0^1\}$. By counting intersections in the corresponding strand diagrams, one sees that the $(M,A)$ bigradings of these generators are $(-1,-1/2)$, $(-1,-1)$, $(0,-1/2)$, $(0,-1/2)$, $(0,0)$, $(-1,-1/2)$, respectively. Label the six empty rectangular regions in the diagram by $A, B, C, D, E, F$, as marked in Figure \ref{fig:cap_hd}, and label the region containing the $O$ by $G$.

Of the six internal rectangles, $A$, $B$, $B\cup G$, $A\cup B$, $G$, and $A\cup B\cup G$, only the first three connect pairs of generators. Of the sets of rectangles intersecting the boundary of the diagram, only the sets consisting of an individual rectangle connect pairs of generators. One can just enumerate all the maps. We provide the resulting type $A$ structure $\ctm(\T_1)$ below. An arrow pointing from  generator $\x$ to a generator $\y$ marked with an algebra element $a$ means that $m_2(\x, a) = \y$. An arrow from $\x$ to $\y$ that is unmarked or marked with $U_1$ means  that $m_1(\x)=\y$ or $m_1(\x) = U_1\y$, respectively. We also provide the respective rectangle for each arrow.

\[
  \mathcenter{
    \begin{tikzpicture}
      \node at (0,0) (a3) {$x_2$};
      \node at (4.8,2.5) (b1) {$y_0$};
      \node at (0,5) (a2) {$x_1$};
      \node at (8,0) (b3) {$y_2$};
      \node at (3.2,2.5) (a1) {$x_0$};
      \node at (8,5) (b2) {$y_1$};
      \draw[algarrow, bend right=20] (b1) to node[pos=0.4,left] {$a_{02}$} node[pos=0.3,right]  {\small{\textcolor{blue}{$C\cup D$}}}  (b3);
      \draw[algarrow, bend left=20] (b1) to node[pos=0.4,left] {$a_{01}$} node[pos=0.3,right] {\small{\textcolor{blue}{$C$}}} (b2);
      \draw[algarrow, bend right=20] (a3) to node[pos=0.6,right] {$a_{20}$} node[pos=0.7,left]  {\small{\textcolor{blue}{$E\cup F$}}}  (a1);
      \draw[algarrow, bend left=20] (a2) to node[pos=0.6,right] {$a_{10}$} node[pos=0.7,left]  {\small{\textcolor{blue}{$E$}}}  (a1);
      \draw[algarrow] (b1) to node[pos=0.5,above]  {\small{\textcolor{blue}{$A$}}}  (a1);
      \draw[algarrow] (a3) to node[pos=0.5,above]  {\small{\textcolor{blue}{$B$}}}  (b3);
      \draw[algarrow] (a2) to node[pos=0.5, above] {$U_1$} node[pos=0.5,below]  {\small{\textcolor{blue}{$G$}}}  (b2);
      \draw[algarrow, bend left=55] (a3) to node[pos=0.5,left] {$a_{21}$} node[pos=0.5,right]  {\small{\textcolor{blue}{$F$}}} (a2);
      \draw[algarrow, bend left=25] (a2) to node[pos=0.5,right] {$U_1 a_{12}$} node[pos=0.5,left] {\small{\textcolor{blue}{$D\cup G$}}} (a3);
      \draw[algarrow, bend right=55] (b3) to node[pos=0.5,right] {$U_1 a_{21}$} node[pos=0.5,left]  {\small{\textcolor{blue}{$G\cup F$}}}  (b2);
            \draw[algarrow, bend right=25] (b2) to node[pos=0.5,left] {$a_{12}$} node[pos=0.5,right] {\small{\textcolor{blue}{$D$}}} (b3);
    \end{tikzpicture}
  }
\]

Using the standard cancelation algorithm for type $A$ structures, we can cancel $y_1$ and $x_1$ to obtain the homotopy equivalent structure
\[
  \mathcenter{
    \begin{tikzpicture}
      \node at (0,0) (a3) {$x_2'$};
      \node at (0,3.5) (a2) {$x_1'$};
      \node at (6,0) (b3) {$y_2'$};
      \node at (6,3.5) (b2) {$y_1'$};
      \draw[algarrow] (a3) to  (b3);
            \draw[algarrow] (a3) to node[pos=0.2, right] {$a_{20}, a_{01}$} (b2);
               \draw[algarrow] (a2) to node[pos=0.2, right] {$a_{10}, a_{02}$} (b3);
        \draw[algarrow, bend right=20] (a3) to node[pos=0.5, below] {$a_{20}, a_{02}$} (b3);
      \draw[algarrow] (a2) to node[pos=0.5, below] {$U_1$}   (b2);
            \draw[algarrow, bend left=20] (a2) to node[pos=0.5, above] {$a_{10}, a_{01}$}   (b2);
      \draw[algarrow, bend left=30] (a3) to node[pos=0.5,left] {$a_{21}$} (a2);
      \draw[algarrow] (a2) to node[pos=0.5,right] {$U_1 a_{12}$}  (a3);
      \draw[algarrow, bend right=30] (b3) to node[pos=0.5,right] {$U_1 a_{21}$}   (b2);
            \draw[algarrow] (b2) to node[pos=0.5,left] {$a_{12}$}  (b3);
    \end{tikzpicture}
  }
\]

Then we can cancel $x_2'$ and $y_2'$ to obtain the structure $\mathcal M$ below:
\[
  \mathcenter{
    \begin{tikzpicture}
      \node at (0,0) (a2) {$x_1''$};
      \node at (4.3,0) (b2) {$y_1''$};
      \draw[algarrow, bend right=15] (a2) to node[pos=0.5, below] {$U_1$}   (b2);
      \draw[algarrow, bend left=15] (a2) to node[pos=0.5, above] {$a_{10}, a_{01}$}   (b2);
      \draw[algarrow, bend left = 55] (a2) to node[pos=0.5, above] {$a_{10}, a_{02}, \vec{a}, \ldots, \vec{a}, a_{20}, a_{01}$}   (b2);
      \draw[algarrow, bend left = 55] (b2) to node[pos=0.5, below] {$a_{12}, \vec{a}, \ldots, \vec{a}, a_{21}$}   (a2);
      \draw[algarrow] (a2) to [out=150, in=210, looseness=9] node[pos=0.5, left] {$a_{10}, a_{02}, \vec{a}, \ldots, \vec{a}, a_{21}$}   (a2);
      \draw[algarrow] (b2) to [out=30, in=330, looseness=9] node[pos=0.5, right] {$a_{12}, \vec{a}, \ldots, \vec{a}, a_{20}, a_{01}$}   (b2);
    \end{tikzpicture}
  }
\]
Here, $\vec{a}$ is the sequence $a_{20}, a_{02}$, and it may be repeated $i\geq 0$ times. For example, the arrow from $x_1''$ to $y_1''$ marked with $a_{12}, \vec{a}, \ldots, \vec{a}, a_{21}$ means that 
\begin{align*}
m_3(x_1'', a_{12}, a_{21}) &= y_1''\\
m_5(x_1'', a_{12}, a_{20}, a_{02}, a_{21}) &= y_1''\\
m_7(x_1'', a_{12}, a_{20}, a_{02}, a_{20}, a_{02}, a_{21}) &= y_1''\\
& \vdots
\end{align*}

\begin{figure}[h]
\centering
  \labellist
  \pinlabel \textcolor{gray}{$J$} at 21 28
  \pinlabel \textcolor{gray}{$K$} at 21 41
  \pinlabel \textcolor{gray}{$N$} at 82 39
  \pinlabel \textcolor{gray}{$L$} at 45 33
  \pinlabel \textcolor{gray}{$M$} at 43 42
  \endlabellist
    \includegraphics[scale=1.25]{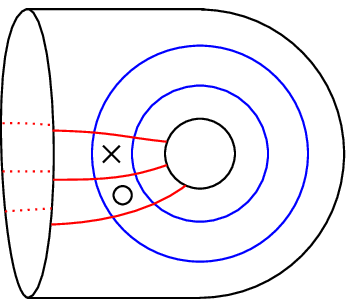}
       \caption{The Heegard diagram for the cap $\T_2$ with $-\bdy^L\T_2 = (-,+)$.}\label{fig:cup_hd}
\end{figure}

Similarly, one can enumerate all rectangles counted in the type $D$ structure maps for the  cap $\T_2$. The Heegaard diagram is displayed in Figure \ref{fig:cup_hd}.  It has six generators: $z_{ij} = \{\alpha_i^0\cap \beta_1^1, \alpha_j^0\cap \beta_2^1\}$, for $\{i,j\}\subset \{0,1,2\}$. The $(M,A)$ bigradings of $z_{01}$, $z_{10}$, $z_{02}$, $z_{20}$, $z_{12}$, and $z_{21}$ are $(-1,-1)$, $(0,0)$, $(-1, -1/2)$, $(0,-1/2)$, $(0,0)$, and $(-1,-1)$, respectively. Label the five empty rectangular regions in the diagram by $J, \ldots, N$, as marked in Figure \ref{fig:cup_hd}, and label the region containing the $O$ by $P$. 

By enumerating the rectangles connecting pairs of generators, one can compute the type $D$ structure  $\ctm(\T_2)$ displayed below. We use the earlier notation for the algebra generators, and we let $U_2$ be the variable corresponding to the $O$ in the Heegaard diagram. If there are $t$ arrows starting at a generator  $\x$ and ending at generators $\y_1, \ldots, \y_t$, marked with algebra elements $a_1, \ldots, a_t$, that means that $\delta^1(\x) = a_1\otimes \y_1+\cdots + a_t\otimes \y_t$.

\[
  \mathcenter{
    \begin{tikzpicture}
      \node at (0,0) (a3) {$z_{21}$};
      \node at (2.8,2.5) (b1) {$z_{01}$};
      \node at (0,5) (a2) {$z_{20}$};
      \node at (8,0) (b3) {$z_{12}$};
      \node at (5.2,2.5) (a1) {$z_{10}$};
      \node at (8,5) (b2) {$z_{02}$};
      \draw[algarrow, bend left=10] (b1) to node[pos=0.5,above] {$a_{21}$} node[pos=0.6,below] {\small{\textcolor{blue}{$M$}}} (b2);
            \draw[algarrow, bend right=15] (a3) to node[pos=0.6,left] {$a_{02}$} node[pos=0.6,right]  {\small{\textcolor{blue}{$J\cup K$}}}  (b1);
      \draw[algarrow, bend left=10] (a2) to node[pos=0.5,above] {$a_{12}$} node[pos=0.4,below]  {\small{\textcolor{blue}{$K$}}}  (a1);
      \draw[algarrow] (b1) to node[pos=0.5,below]  {$U_2$} node[pos=0.5,above]  {\small{\textcolor{blue}{$P$}}}  (a1);
      \draw[algarrow] (a3) to  node[pos=0.5,below] {$1+U_2$} node[pos=0.5,above]  {\small{\textcolor{blue}{$K\cup M, N\cup P$}}}  (b3);
      \draw[algarrow] (a2) to node[pos=0.5, above] {$1$} node[pos=0.5,below]  {\small{\textcolor{blue}{$N$}}}  (b2);
      \draw[algarrow, bend left=50] (a3) to node[pos=0.5,left] {$U_2 a_{01}$} node[pos=0.5,right]  {\small{\textcolor{blue}{$J\cup P$}}} (a2);
      \draw[algarrow, bend left=25] (a2) to node[pos=0.5,right] {$a_{10}$} node[pos=0.5,left] {\small{\textcolor{blue}{$L$}}} (a3);
      \draw[algarrow, bend right=50] (b3) to node[pos=0.5,right] {$a_{01}$} node[pos=0.5,left]  {\small{\textcolor{blue}{$J$}}}  (b2);
            \draw[algarrow, bend right=25] (b2) to node[pos=0.5,left] {$U_2 a_{10}$} node[pos=0.5,right] {\small{\textcolor{blue}{$L\cup P$}}} (b3);
            \draw[algarrow, bend right=15] (a1) to node[pos=0.4,right] {$a_{20}$} node[pos=0.4,left]  {\small{\textcolor{blue}{$L\cup M$}}}  (b3);
    \end{tikzpicture}
  }
\]

Since $\ctm(\T_2)$ is bounded, we can take the box tensor product of any right type $A$ structure over $\mathcal A(-,+)$ with it.
The chain complex $\mathcal M\boxtimes\ctm(\T_2)$ is generated by $x_1''\boxtimes z_{02}$, $x_1''\boxtimes z_{20}$, $y_1''\boxtimes z_{02}$, and $y_1''\boxtimes z_{20}$, in $(M, A)$ bigradings $(-2,-3/2)$, $(-1,-3/2)$, $(-1,-1/2)$, and $(0,-1/2)$, respectively. By pairing type $A$ and type $D$ maps, we see that the differential is given by 
\begin{align*}
d(x_1''\boxtimes z_{20}) &= (U_1+U_2)y_1''\boxtimes z_{20},\\
d(x_1''\boxtimes z_{02}) &= (U_1+U_2)y_1''\boxtimes z_{02},\\
d(y_1''\boxtimes z_{20}) &= 0,\\
d(y_1''\boxtimes z_{02}) &= 0.
\end{align*}
As a complex over $\Ft[U_1]$, this is homotopy equivalent to $\Ft[U_1]\otimes (\Ft\oplus\Ft)$ , generated by $y_1''\boxtimes z_{20}$ and $y_1''\boxtimes z_{02}$. After shifting bigradings by $(l/2, l/2)$, this agrees with $\CFKm(\mathcal U) \otimes (\Ft[-1/2]\oplus\Ft[1/2])$.


\end{document}